\documentclass{aimscleaned}

\usepackage{graphics} 
\usepackage{epsfig} 
\usepackage{graphicx}  \usepackage{epstopdf}
\usepackage[colorlinks=true]{hyperref}
\hypersetup{urlcolor=blue, citecolor=red}

\def\currentvolume{X}

    \def\ppages{X--XX}
     \def\DOI{10.3934/xx.xx.xx.xx}

\newtheorem{theorem}{Theorem}[section]

\theoremstyle{definition}

\newtheorem{remark}{Remark}

\usepackage{color}
\usepackage[percent]{overpic}

\newcommand{\R}{\mathbb{R}}

\newcommand{\Dt}{\Delta t}
\newcommand{\yii}{y_{i+1}}
\newcommand{\yi}{y_i}
\newcommand{\wmax}{v_{\textup{max}}}
\newcommand{\vmax}{v_{\textup{max}}}
\newcommand{\JJ}{\textsc{j}}
\newcommand{\RR}{\textsc{r}}
\newcommand{\next}{\textsc{next}}
\newcommand{\yip}{y_{i,p}}
\newcommand{\ykt}{y_k(t)}
\newcommand{\ykkt}{y_{k+1}(t)}
\newcommand{\dikt}{\delta_{i_k}(t)}
\newcommand{\dit}{\delta_{i}(t)}
\newcommand{\Dkt}{\Delta_{k}(t)}
\newcommand{\Dkit}{\Delta_{k_i}(t)}
\newcommand{\baryii}{\bar y_{i+1}}
\newcommand{\baryi}{\bar y_i}
\newcommand{\bardi}{\bar \delta_i}
\newcommand{\dykt}{\dot y_k(t)}
\newcommand{\Cphi}{C_\varphi}
\newcommand{\Ouno}{\mathcal{O}(1)}
\newcommand{\ykit}{y_{k_i}(t)}
\newcommand{\yiit}{y_{i+1}(t)}
\newcommand{\yit}{y_{i}(t)}

\title[On the micro-to-macro limit for traffic flow on networks] 
{On the micro-to-macro limit for first-order traffic flow models on networks}

\author[E. Cristiani and S. Sahu]{}


 \email{e.cristiani@iac.cnr.it}
 \email{sahu@mat.uniroma1.it}

\begin{document}
\maketitle

\centerline{\scshape Emiliano Cristiani}
{\footnotesize
 \centerline{Istituto per le Applicazioni del Calcolo ``M. Picone''}
   \centerline{Consiglio Nazionale delle Ricerche}
   \centerline{Rome, Italy}
} 

\medskip

\centerline{\scshape Smita Sahu}
{\footnotesize
 \centerline{Dipartimento di Matematica ``G. Castelnuovo''}
   \centerline{Sapienza, Universit\`a di Roma}
   \centerline{Rome, Italy}
} 
\bigskip


\begin{abstract}
Connections between microscopic follow-the-leader and macroscopic fluid-dynamics traffic flow models are already well understood in the case of vehicles moving on a single road. Analogous connections in the case of road networks are instead lacking. This is probably due to the fact that macroscopic traffic models on networks are in general ill-posed, since the conservation of the mass is not sufficient alone to characterize a unique solution at junctions. This ambiguity makes more difficult to find the right limit of the microscopic model, which, in turn, can be defined in different ways near the junctions.
In this paper we show that a natural extension of the first-order follow-the-leader model on networks corresponds, as the number of vehicles tends to infinity, to the LWR-based multi-path model introduced in \cite{bretti2014DCDS-S, briani2014NHM}.
\end{abstract}

\section{Introduction}\label{sec:intro}
Traffic flow can be described at different scales, depending on the level of details one wants to observe. Typically, three scales of observation can be adopted: \textit{microscopic} (single vehicles are tracked), \textit{mesoscopic} (averaged quantities like density and mean velocity are tracked but car-to-car interactions are not completely lost) and \textit{macroscopic} (only averaged quantities are observed). 

Limiting our attention to differential models, the most famous microscopic models are those of \textit{follow-the-leader} type, also known as \textit{car-following} models \cite{helbing2001RMP, pipes1953JAP}. In such a models, the dynamics of each vehicle depend on the vehicle in front of it, so that, in a cascade, the whole traffic flow is determined by the dynamics of the very first vehicle (the \textit{leader}).
The oldest macroscopic model is instead the LWR model \cite{LW,R}, which dates back to 1955 and is inspired by the laws of fluid-dynamics.

Generally speaking, microscopic models are considered more justifiable because the behavior of every single vehicle can be described with high precision and it is immediately clear which kind of interactions are considered. On the contrary, macroscopic models are based on assumptions which are hardly correct or at least verifiable (among others, the actual validity of the continuum hypothesis and the negligibleness of car-to-car interactions). As a consequence, it is often desirable establishing a connection between microscopic and macroscopic models in order to justify and validate the latter on the basis of the verifiable modeling assumptions of the former.

Connections between microscopic follow-the-leader and macroscopic fluid-dynamics traffic flow models are already well understood in the case of vehicles moving on a \emph{single road}. Analogous connections in the case of \emph{road networks} are instead lacking. This is probably due to the fact that macroscopic traffic models on networks are in general ill-posed, since the conservation of the mass is not sufficient alone to characterize a unique solution at junctions. This ambiguity makes more difficult to find the right limit of the microscopic model, which, in turn, can be defined in different ways near the junctions.

\emph{Goal.}
In this paper we propose a very natural extension of a first-order follow-the-leader model on road networks and then we prove that its solution tends to the solution of the LWR-based multi-path model introduced in \cite{bretti2014DCDS-S, briani2014NHM} in the limit, i.e.\ as the number of vehicles tends to infinity while their total length is kept constant. The limit is proved extending to networks the results already existing for a single road, and it is then confirmed by numerical experiments. It is useful to recall that the multi-path method is able to select automatically an admissible solution at junction, thus resolving ill-posedness issues. We also recall that the solution selected by the multi-path method does not match the one obtained by maximizing the flux at junction. Therefore, the connection with the microscopic model promotes the solution computed by the multi-path method as ``more natural'', while the one computed by maximizing the flux should be seen as the ``most desirable'', to be achieved by means of \emph{ad hoc} traffic regulations.

\emph{Relevant literature.} 
The literature about microscopic and macroscopic traffic flow models is huge and a detailed review is out of the scope of the paper. For a quick introduction to the field we suggest the book \cite{habermanbook} and the surveys \cite{bellomo2011SR, helbing2001RMP}. 
Regarding first-order models on a single unidirectional road, the micro-to-macro limit was already deeply investigated by means of different techniques:
the papers \cite{colombo2014RSMUP, rossi2014DCDS-S} use standard techniques coming from the theory of conservation laws;   
the paper \cite{difrancesco2015ARMA} instead proves the limit relying also on measure theory. The microscopic solution is interpreted as an empirical measure which is proven to converge to the entropy solution of the macroscopic model in the 1-Wesserstein topology;  
finally, the papers \cite{costeseque2011thesis, forcadel2014hal} attack the problem exploiting the link between conservation laws and Hamilton-Jacobi equations. 

Micro-to-macro limit for second-order models was instead investigated in \cite{aw2002SIAP, greenberg2001SIAP}, where the Aw-Rascle model is derived as the limit of a second-order follow-the-leader model.

Macroscopic-only traffic models on networks were deeply investigated starting from \cite{holden1995SIMA}. 
A complete introduction can be found in the book \cite{piccolibook}, which discusses several methods to characterize a unique solution at junctions. Let us also mention the source-destination model introduced in \cite{garavello2005CMS} (see also \cite{herty2008CMS}) and the buffer models \cite{garavello2012DCDS-A, garavello2013bookchapt, herty2009NHM}. 
Recently, the LWR-based multi-path model on networks was introduced in the paper \cite{bretti2014DCDS-S}, together with a Godunov-based numerical scheme to solve the associated system of conservation laws with discontinuous flux. The relationship between the multi-path model and more standard methods (like, e.g., the one proposed in \cite{coclite2005SIMA} based on the maximization of the flux at junction) was investigated in \cite{briani2014NHM}.

To our knowledge, there is no systematic theory about the extension of the follow-the-leader models on networks. It is plain that at the microscopic level one can easily reach a high level of detail, including junctions with spatial extension (non point), multi-lane roads, multi-class vehicles, traffic lights and priorities. Several highly sophisticated simulators are available since many years (free and commercial), see, e.g., \cite{felleVISSIM2010, shen2012GM} and references therein to have an idea of the models and methods commonly used. Nevertheless, it is unclear which average flux is actually observed at junction by the many-particle limit of any car-following model. 

Beside this, let us also mention that the relationship between microscopic and macroscopic models was exploited to create hybrid models, see, e.g., \cite{moutari2007SIAP}. In such a models the averaged quantities are observed where a detailed description is not needed (e.g., far from the junctions) while microscopic dynamics are considered elsewhere. However, this approach gives no clue about the macroscopic behavior of the microscopic model at junctions.

\emph{Paper organization.} 
In section~\ref{sec:background} we introduce the follow-the-leader model and the LWR model on a single road. We also recall the existing results about the micro-to-macro limit on a single road following \cite{colombo2014RSMUP, difrancesco2015ARMA, rossi2014DCDS-S}. 
In section \ref{sec:FTLnetw} we extend the follow-the-leader model to networks, and in section \ref{sec:limit}, which is the core of the paper, we show the relationship between the follow-the-leader model on networks and the LWR-based multi-path model. 
In section \ref{sec:numerics} we present fully discrete algorithms for the numerical solution to the equations related to the models previously discussed and finally in section \ref{sec:tests} we confirm our findings by means of some numerical tests.

\section{Background and previous results}\label{sec:background}
Let us describe the first-order follow-the-leader model we will use as main ingredient in the rest of the paper. We assume here that vehicles move on a single infinite road with a single lane. Vehicles are initially located one after the other and cannot overtake each other. We denote by $n\in\mathbb N$ the number of vehicles, by $\ell_n>0$ the length of the single vehicles, and by $\mathcal L$ their total length. We have
\begin{equation}
\ell_n=\frac{\mathcal L}{n}.
\end{equation} 
This relation is crucial for the micro-to-macro limit since it translates the fact that the total length does not change when the number of vehicles tends to infinity, because the cars ``shrink'' accordingly.

We denote by $y_i(t)$ the position of the $i$-th car and we assume that at the initial time $t=0$ cars are labeled in order, i.e.\ $y_1(0)<y_2(0)<\ldots <y_n(0)$. This guarantees that the $(i+1)$-th car is just in front of the $i$-th one. Moreover, it is assumed that at $t=0$ cars do not overlap, i.e.\ $\yii(0) - \yi(0)\geq\ell_n$, $i=1,\ldots,n-1$. We are now ready to introduce the model, described by the following system of ODEs
\begin{equation}\label{FTL_1road}
\left\{
\begin{array}{ll}
\dot y_i(t)=w(\delta_i(t)), & i=1,\ldots,n-1 \\
\dot y_n(t)=\wmax,
\end{array}
\right.
\end{equation}
where 
$$
\delta_i(t):=\yii(t)-\yi(t)
$$
and $w$ is such that 
\begin{equation}\label{link_v_w}
w:[\ell_n,+\infty)\to[0,\wmax], \qquad w(\delta):=v\left(\frac{\ell_n}{\delta}\right),
\end{equation}
with $v\in C^1([0,1];[0,\wmax])$ any function such that  
\begin{equation}\label{properties_of_v}
v'(r) < 0, \qquad v(0)=\wmax, \qquad v(1)=0.
\end{equation}
The $n$-th car is the leader and it is assumed to move at maximal velocity $\wmax>0$. 
Note that the properties introduced above guarantee that the cars do not overlap at any later time $t>0$, see \cite[Lemma 1]{difrancesco2015ARMA}.

\medskip

The macroscopic limit of the previous model is given by the well known LWR model, which describes the evolution of the average (normalized) density of vehicles $\rho(t,x):([0,+\infty)\times\R)\to [0,1]$ by means of the following conservation law
\begin{equation}\label{LWR_1road}
\partial_t\rho + \partial_x(\rho v(\rho))=0, \qquad  t>0, \quad x\in\R.
\end{equation}
Relationship (\ref{link_v_w}) makes the link between the two models.

In order to recall precisely the results about the correspondence between the two models, we need to introduce first the natural spaces for the macroscopic density $\rho$ and for the vectors
of vehicles' positions $\mathbf y=(y_1,\ldots,y_n)$, at any fixed time:
$$
R:=\left\{r\in L^1(\R;[0,1])~:~\int_\R r(x)dx=\mathcal L \text{ and supp($r$) is compact}\right\}
$$
and
$$
Y_n:=\Big\{\mathbf{y}\in\R^n~:~\yii-\yi\geq\ell_n, \quad \forall i=1,\ldots,n-1\Big\}.
$$

We also introduce the  operators $E_n:R\to Y_n$ and $C_n:Y_n\to R$, defined respectively as
\begin{equation}\label{def:En}
E_n[r(\cdot)]:=\mathbf{y}=
\left\{
\begin{array}{l}
y_n=\max(\text{supp}(r)), \\
y_i=\max\left\{z\in\R~:~\displaystyle\int_z^{\yii}r(x)dx=\ell_n\right\}, \quad i=n-1,\ldots,2,1
\end{array}
\right.
\end{equation}
and
\begin{equation}\label{def:Cn}
C_n[\mathbf{y}]:=r_n=\sum_{i=1}^{n-1}\frac{\ell_n}{\delta_i}\chi_{[\yi,\yii)},
\end{equation}
where $x\to \chi_I(x)$ is the indicator function of any subset $I\subset\R$. 
The discretization operator $E_n$ acts on a macroscopic density $\rho(t,\cdot)$, providing a vector of positions $\mathbf{y}(t)=E_n[\rho(t,\cdot)]$ whose components partition the support of the density into
segments on which $\rho$ has fixed integral $\ell_n$. 
On the contrary, the operator $C_n$ antidiscretizes a microscopic vector of positions $\mathbf{y}$, giving a \textit{piecewise constant} density $\rho_n(t,\cdot)=C_n[\mathbf{y}(t)]$.

We are now ready to state the main result about the convergence of the microscopic model (\ref{FTL_1road}) to the macroscopic one (\ref{LWR_1road}). The proof can be found in \cite{difrancesco2015ARMA} (see also \cite{colombo2014RSMUP}).

\begin{theorem}\label{teo:rossi}
Let (\ref{link_v_w}) and (\ref{properties_of_v}) hold. Choose $\bar\rho=\rho(0,\cdot)\in R~\cap~BV(\R;[0,1])$ and $\bar{\mathbf y}=\mathbf y(0)=E_n[\bar\rho]$. 
Let $\mathbf y(\cdot)$ be the solution of (\ref{FTL_1road}) with initial condition $\bar{\mathbf y}$.
Define $\rho_n(t,\cdot)=C_n[\mathbf y(t)]$. Then $\rho_n$ converges almost everywhere and in $ L^1_{\textup{loc}}([0,+\infty)\times\R)$ to the unique entropy solution $\rho$ to the problem (\ref{LWR_1road}) with initial condition $\rho(0,\cdot)=\bar\rho$.
\end{theorem}
%
%
%
%
\section{Follow-the-leader model on networks}\label{sec:FTLnetw}
In this section we extend the follow-the-leader model described in section \ref{sec:background} to a road network. We define a road network $\mathcal G$ as a direct graph with $N_\RR$ arcs (roads) and $N_\JJ$ nodes (junctions). 
We assume that for each junction $\textsc{j}=1,\ldots,N_\JJ$ , there exist disjoint subsets $\textsc{inc}(\textsc{j}), \textsc{out}(\textsc{j}) \subseteq \{1,\ldots,N_\RR\}$, representing, respectively, the incoming roads to $\textsc{j}$ and the outgoing roads from $\textsc{j}$. Among junctions, we distinguish two particular subsets consisting of origins $\mathcal O$, which are the junctions $\textsc{j}$ such that $\textsc{inc}(\textsc{j})=\emptyset$, and destinations $\mathcal D$, which are the junctions $\textsc{j}$ such that $\textsc{out}(\textsc{j})=\emptyset$.
Finally, we denote by $L^\RR$ the length of the road $\RR$ for any $\RR=1,\ldots,N_\RR$. At initial time, $n$ vehicles are located anywhere in the network. 

\subsection{A natural extension}\label{sec:FTLnetw_natural}
A natural extension of the follow-the-leader model is derived as follows.
We label the $n$ vehicles by an index $i=1,\ldots,n$ and we denote by $\RR(i,t)\in\{1,\ldots,N_\RR\}$ the road that the vehicle $i$ is traveling at time $t$. We also denote by $y_i(t)$ the position of the vehicle $i$ at time $t$, defined as the distance from the origin of the road $\RR(i,t)$.

There are three main differences with respect to the model on a single road, which require some modifications in the definition of $\delta_i$ and $w$.

\textbf{A}. 
The concept of ``ahead'' must be redefined because at junctions it is unclear which is the space in front of a car, unless a preference among the outgoing roads is assigned. We assume that each vehicle has a sequence of roads (i.e., a path) assigned at the initial time to travel along, so that it is always clear which is the space in front of the vehicle. We also assume that drivers can actually see the space ahead, so that they can always evaluate the distance between them and the car in front, even if the latter is in the next road. Given that, \emph{a vehicle is a leader if it is the first vehicle on its road and there is no other vehicles in the next roads along its path}. It is important to note that now more than one leader can be present on the network at the same time, and vehicles can gain or loose the leadership. 

\textbf{B}. 
The vehicle $i+1$ is no longer necessarily in front of vehicle $i$. This force us to introduce a new notation to refer to the vehicle in front. If a vehicle is not a leader, we denote by $\next(i)$ the label of the car just in front of car $i$. If it is a leader, we set $\next(i)=\emptyset$.

The distance between a (nonleader) vehicle and the vehicle in front is now computed as
\begin{equation}\label{def:deltai}
\hat\delta_i(t):=\left\{
\begin{array}{ll}
y_{\next(i)}(t)-\yi(t), & \text{if } \RR(i,t)=\RR(\next(i),t)  \\ [2mm]
(L^{\RR(i,t)}-\yi(t))+\sum_{\textsc{s}} L^\textsc{s} + y_{\next(i)}(t), & \text{if } \RR(i,t)\neq\RR(\next(i),t)
\end{array}
\right.
\end{equation}
where the summation with respect to \textsc{s} is done over all the empty roads between $\RR(i,t)$ and $\RR(\next(i),t)$ (if any) along the sequence of roads followed by vehicle $i$. 
It is important to note that $\hat\delta_i(t)$ is in general discontinuous because at any time other cars can join or leave the $i$-th car's path, thus reducing or increasing abruptly the distance $\hat\delta_i$ (or even changing the leader/follower status of the $i$-th car).

\begin{figure}[h!]
\begin{center}
\includegraphics[width=0.49\textwidth]{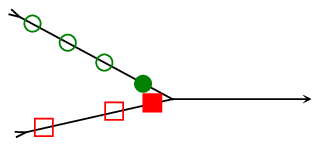}
\includegraphics[width=0.49\textwidth]{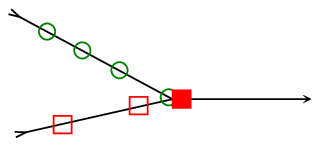}
\end{center}
\caption{Overlapping issues across the junction. Leaders are denoted by filled square/circle.}
\label{fig:overlapping}
\end{figure}

\textbf{C}. 
The non-overlapping condition $\hat\delta_i\geq\ell_n$ is no longer guaranteed, even if it holds at the initial time $t=0$. To see this, let us consider the simple case of a merge, see Fig.\ \ref{fig:overlapping}. Assume that the outgoing road is empty and some cars are traveling along both the incoming roads (Fig.\ \ref{fig:overlapping}-left). Assume also that the two leader cars are going to reach the junction nearly or exactly at the same moment (note that two cars are leader at this time). As soon as one of the two cars goes through the junction, the other one becomes a follower and its distance $\hat\delta_i$, defined by (\ref{def:deltai}), could be less than $\ell_n$ or even 0 (Fig.\ \ref{fig:overlapping}-right).
However, two important observations can be made:
\begin{itemize}
\item[C1.] The number of overlapping cars at a junction is bounded by the number of incoming roads of that junction. In other words, cars do not stack one over the other without limits, since the newly coming vehicles perceive the presence of the overlapping vehicles in front and behave normally, without getting too close to them. As a consequence, the total length of overlapping vehicles will tend to zero as $n\to +\infty$.
\item[C2.] The space interval in which overlap can occur is $(L^\RR-\ell_n,L^\RR]$, which shrinks to $\emptyset$ when $n\to+\infty$.
\end{itemize}
From now on we allow the above described ``mild'' overlap to happen and we handle it extending the function $w$ defined in (\ref{link_v_w}) by means of a new function $w^*$ defined as follows:
\begin{equation}\label{def:w*}
w^*:[0,+\infty)\to[0,\wmax], \qquad 
w^*(\delta)=\left\{
\begin{array}{ll}
w(\delta), & \text{if } \delta\geq\ell_n\\
0, & \text{if } \delta\leq\ell_n 
\end{array}
\right.
\end{equation}
i.e.\ we assume that (smashed) vehicles with $\hat\delta_i<\ell_n$ stops completely until the vehicle in front leaves a space $>\ell_n$, then they re-start moving normally. If more than one vehicle is found to be \emph{exactly} in the same place, the one with the largest label is taken as vehicle ``in front'' and the others ``behind''.

\medskip

We are now ready to introduce the follow-the-leader model on networks. For any $i=1,\ldots,n$ we write
\begin{equation}\label{FTL_netw}
\left\{
\begin{array}{ll}
\dot y_i(t)=w^*(\hat\delta_i(t)), & \text{if } \next(i)\neq\emptyset \\
\dot y_i(t)=\wmax, & \text{if } \next(i)=\emptyset
\end{array}
\right.
\end{equation}
with suitable admissible initial conditions $y_i(0)=\bar y_i$.
In addition, when $y_i(t)=L^{\RR(i,t)}$, the road $\RR(i,t)$ must be updated according to the $i$-th vehicle's path and $y_i(t)$ must be reset to 0 (similarly to a new initial condition).

\begin{remark}
The extension of the follow-the-leader model described above defines implicitly a certain behavior of vehicles at junctions. Since vehicles cannot overtake each other, they basically adopt a FIFO behavior (see \cite[Appendix B]{herty2003SISC} for a general discussion in the context of traffic flow). Consider for example the case of a diverge (see Fig.\ \ref{fig:networks_for_tests}(top-right)) and assume that one of the two outgoing roads is fully congested. If the first car on the incoming road wants to go to the congested road, it will stop. After that, all the following cars will stop too, even if they want to turn to the other outgoing road. The FIFO behavior will be kept in the macroscopic limit.
\end{remark}
\begin{remark}
Beside the overlap issue, the model exhibits other unrealistic behaviors. Let us consider the junction depicted in Fig.\ \ref{fig:networks_for_tests}(bottom). If two cars arrive at the junction from distinct roads exactly at the same time and then they want to go to distinct roads, they will not be affected at all by each other, and, in particular, they will not slow down. This odd behavior can be easily avoided by considering enlarged (non point) junctions or allowing cars which approach a junction to see along \textit{all} the incoming roads, but then it can be more difficult to prove the relationship between the microscopic model and the corresponding macroscopic model.
\end{remark}
\subsection{The model reformulated}\label{sec:FTLnewtpaths}
\begin{figure}[h!]
\begin{center}
\begin{overpic}
[width=0.49\textwidth]{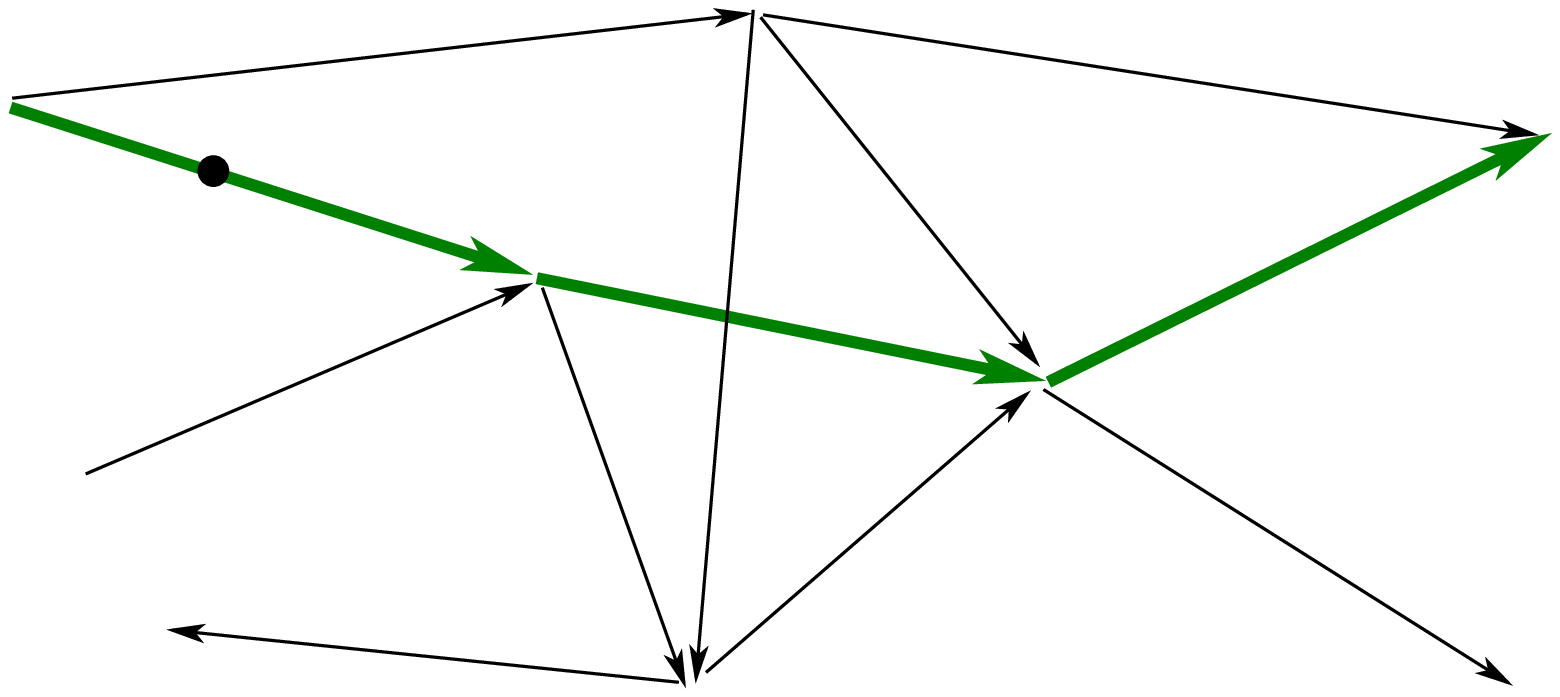}
\put(8,28){\footnotesize $(i,p)$}
\end{overpic}
\begin{overpic}
[width=0.49\textwidth]{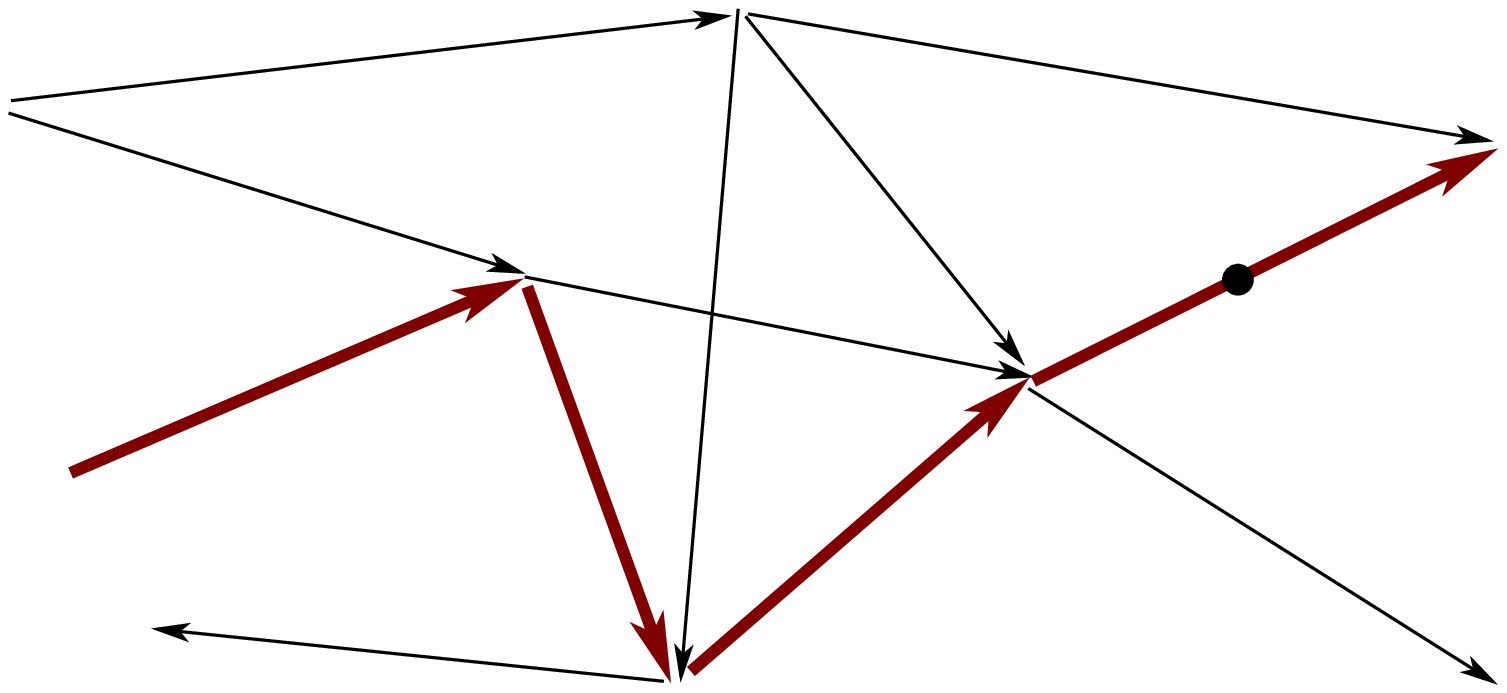}
\put(80,22){\footnotesize $(k,q)$}
\end{overpic}
\end{center}
\caption{A generic network $\mathcal G$. Two possible paths are highlighted. Note that the two paths share an arc of the network.}
\label{fig:paths}
\end{figure}

In order to pass to the limit for $n\to +\infty$, it is convenient to reformulate the model introduced in the previous section. Let us consider all the possible paths on $\mathcal G$ joining all the origin nodes $\mathcal O$ with all the destination nodes $\mathcal D$, see Fig.\ \ref{fig:paths}. Each path is considered as a \emph{single uninterrupted road}, with no junctions (except for its origin and destination). 
Note that paths can share some arcs of the network. 
Let us denote by $P$ the total number of paths and let us divide the $n$ vehicles in $P$ populations, on the basis of the path they are following. 

Let us denote by $n_p$, $p=1,\ldots,P$, the number of vehicles following path $p$ (we have $\sum_{p=1}^{P}n_p=n$), and label \emph{univocally} all vehicles by the multi-index $(i,p)$, $i=1,\ldots,n_p$, \ $p=1,\ldots,P$. 
Let us also denote by $\yip(t)$ the position of the vehicle $i$ of population $p$ at time $t$, defined as the distance from the origin of the path $p$ (not from the origin of the current road). Since paths can overlap, a vehicle $(k,q)$ belonging to a population $q\neq p$ with position $y_{k,q}(t)$ can be found along path $p$ at time $t$. Let us denote by $c_p(k,q,t)$ the distance between the vehicle $(k,q)$ at time $t$ and the origin of path $p$, along path $p$. In other words $c_p$ acts as a change of coordinates between any path $q\neq p$ intersecting path $p$ and path $p$. The function $c_p$ is trivially extended to vehicles already belonging to population $p$ simply setting $c_p(i,p,t)=y_{i,p}(t)$.

Most important, we redefine the concept of distance between two vehicles. 
First, we define $\next(i,p;q)$ as the multi-index of the nearest vehicle belonging to population $q$ in front of the vehicle $(i,p)$ along the path $p$ (setting it to $\emptyset$ if there is no such a vehicle). Then, we define the $p$-\emph{distance} between a (nonleader) vehicle and the vehicle in front as
$$
\delta_i^p(t):=y_{\next(i,p;p)}(t)-\yip(t).
$$
This corresponds to consider only vehicles of population $p$, neglecting the presence of the others. 
We also keep considering the distance between two contiguous vehicles traveling along path $p$, regardless the population they belong to. To this end we define, for any nonleader vehicle of any population,
$$
\Delta_{k,q}^p(t):=c_p(\next(k,q;\bullet),t)-c_p(k,q,t),
$$
where $\next(k,q;\bullet)$ is the vehicle in front of the vehicle $(k,q)$, no matter the population it belongs to. 

In the new formulation the model (\ref{FTL_netw}) becomes,
\begin{equation}\label{FTL_netw_paths}
\left\{
\begin{array}{ll}
\dot y_{i,p}(t)=w^*(\Delta_{i,p}^p(t)), & \text{if } \next(i,p;\bullet)\neq\emptyset \\
\dot y_{i,p}(t)=\wmax, & \text{if } \next(i,p;\bullet)=\emptyset
\end{array}
\right.
\end{equation}
for any $p=1,\ldots,P$ and $i=1,\ldots,n_p$. 
Note that (\ref{FTL_netw_paths}) is a system of $P$ coupled systems of ODEs with discontinuous right-hand side. As already mentioned, the discontinuity comes from the fact that cars following paths other than $p$ can abruptly join or leave path $p$, thus decreasing or increasing the distance $\Delta_{i,p}^p$. However, in the time intervals in which no car crosses junctions the system falls in the standard theory of the ODEs. To deal with the discontinuity, we will consider the integral form of (\ref{FTL_netw_paths}) and we assume that there exists a solution in the sense of Carath\'eodory.

\section{Micro-to-macro limit}\label{sec:limit}
\begin{figure}[t!]
\vskip1cm
\begin{center}
\begin{overpic}
[width=0.75\textwidth]{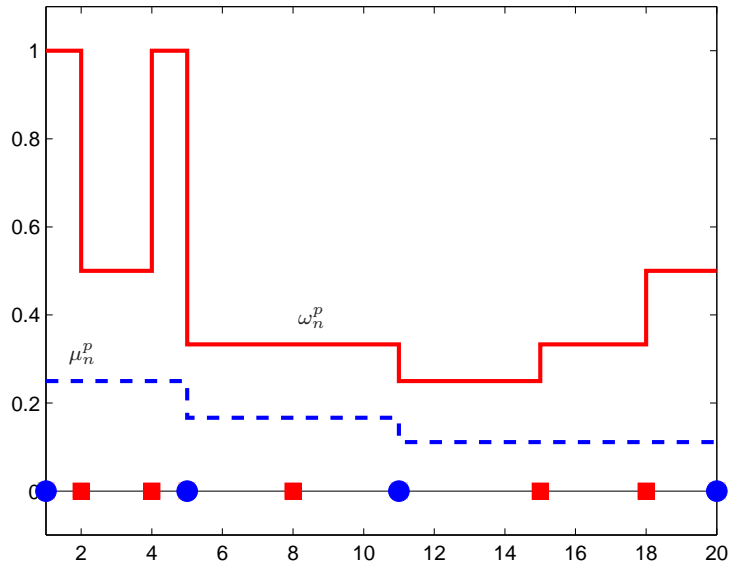}
\put(8,28){\footnotesize $\mu^p_n$} \put(40,33){\footnotesize $\omega^p_n$}
\end{overpic}
\end{center}
\caption{Construction of functions $\mu_n^p$ (dashed blue) and $\omega_n^p$ (solid red). Cars belonging to population $p$ are marked by blue circles, all the others by red squares. Here $\mathcal L=1$, $n=9$, $n_p=4$.}
\label{fig:mu_omega}
\end{figure}

To begin with, we note that condition (\ref{link_v_w}) clarifies the relationship between the average density and the distance between vehicles which is valid in the limit $n\to +\infty$ (see Theorem \ref{teo:rossi}). It states that
\begin{equation}\label{m2M_crucial}
\rho(t,\cdot)\longleftrightarrow\frac{\ell_n}{\delta(t)}.
\end{equation}
In section \ref{sec:FTLnewtpaths} we have introduced $p$-distances which must be now related to the right average densities. To do that, let us use again (with an abuse of notation) the operators $E_n$ and $C_n$ introduced in section \ref{sec:background}, modifying in a obvious manner the sets $R$ (considering the total length $\mathcal L^p$ of cars belonging to population $p$) and $Y_n$ to deal with the new framework. 
Define
$$
\mathbf{y}^p:=(y_{1,p},\ldots,y_{n_p,p}), \quad\text{ for any } p=1,\ldots,P ,
$$
and
$$
\mathbf{y}:=(\mathbf{y}^1,\ldots,\mathbf{y}^{P}).
$$
The new density functions are
$$ 
\mu^p_n(t,\cdot):=C_n[\mathbf{y}^p(t);\mathcal G,\delta^p]= \sum_{\substack{\text{followers $(i,p)$} \\ 
\text{of pop.\ } p}}\frac{\ell_n}{\delta_i^p}\chi_{[y_{i,p}(t),y_{\next(i,p;p)}(t))},
$$
$$
\omega^p_n(t,\cdot):=C_n[\mathbf{y}(t);\mathcal G,\Delta^p]= \sum_{\substack{\text{followers $(k,q)$} \\ 
\text{on path } p \\ \text{(any pop.)}}}\frac{\ell_n}{\Delta_{k,q}^p}\chi_{[c_p(k,q,t),c_p(\next(k,q;\bullet),t))}.
$$
Both functions $\mu^p_n(t,\cdot)$ and $\omega^p_n(t,\cdot)$ are defined on the whole path $p$. The function $\mu^p_n$ corresponds to the density of vehicles belonging to population $p$, while the function $\omega^p_n$ represents the \emph{total} density along path $p$, which depends on all vehicles belonging to population $p$ plus all vehicles following a path $q$ which shares some roads with path $p$, see Fig.\ \ref{fig:mu_omega}.

Defining the limits (if they exist)
$$
\mu^p(t,x):=\lim_{n\to +\infty} \mu^p_n(t,x), \qquad \omega^p(t,x):=\lim_{n\to +\infty} \omega^p_n(t,x),
$$
we get, in analogy to (\ref{m2M_crucial}),
\begin{equation}\label{m2M_new}
\mu^p(t,\cdot)\longleftrightarrow\frac{\ell_n}{\delta^p(t)}, \qquad
\omega^p(t,\cdot)\longleftrightarrow\frac{\ell_n}{\Delta^p(t)}.
\end{equation}

In order to introduce the macroscopic model we will also need the correspondence between the microscopic and the macroscopic velocity. To this end, we define the function $v^*$ such that
\begin{equation}\label{def:v*}
v^*:[0,+\infty)\to[0,\vmax],\qquad 
v^*(r)=
\left\{
\begin{array}{ll}
v(r), & \text{if } r\leq 1 \\
0, & \text{if } r\geq 1.
\end{array}
\right.
\end{equation}
It is easy to see that the relation in (\ref{link_v_w}) still holds, more precisely
$$
w^*(\delta)=v^*\left(\frac{\ell_n}{\delta}\right).
$$

At this point it is crucial to note that a path, by definition, has no junctions (except for its origin and destination) and therefore is indistinguishable from a single road. As a consequence, we can follow in broad terms the results already proved in \cite{colombo2014RSMUP, difrancesco2015ARMA, rossi2014DCDS-S} for the micro-to-macro limit on a single road.

Using (\ref{m2M_new}), and reasoning by analogy with the 1D case, we claim that, for any fixed $p=1,\ldots,P$, the macroscopic equation associated to the model (\ref{FTL_netw_paths}) is
\begin{equation}\label{LWRlimit}
\partial_t\mu^p + \partial_x\left(\mu^p v^*(\omega^p)\right)=0, \qquad t>0, \quad x\in\R.
\end{equation}
Note that in equation (\ref{LWRlimit}) both densities $\mu^p$ and $\omega^p$ appear. This is due to the fact that equation (\ref{FTL_netw_paths}) describes the evolution of the population $p$ only, but the velocity is evaluated considering the distance $\Delta^p$, which is indeed related to the \emph{total} density on the path. Equation (\ref{LWRlimit}) accounts for all populations of cars and then it must be seen as a system of PDEs.

\medskip

To prove the correspondence, let us fictitiously extend each path to $(-\infty,+\infty)$, setting the density to 0 after the leader and before the last follower. In this way we simplify the problem getting rid of the boundary conditions. 

\begin{figure}[h!]
\begin{center}
\begin{overpic}
[width=\textwidth]{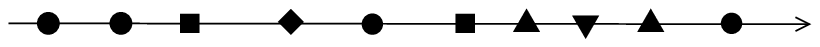}
\put(12.0,5){$i\!-\!1$}\put(45,5){$i$}\put(88,5){$i\!+\!1$}\put(63.5,5){$k_i$}
\end{overpic}\\ \vskip1cm
\begin{overpic}
[width=\textwidth]{./figs/proof}
\put(44.5,5){$i_k$}\put(63.5,5){$k$}\put(69.5,5){$k\!+\!1$}\put(87,5){$i_k\!+\!1$}
\end{overpic}\end{center}
\caption{Path $p$. Definition of $k_i$ (top) and definition of $i_k$ (bottom).}
\label{fig:proof}
\end{figure}

Before entering the proof of the main result, we need to simplify the notations. Fix $p\in\{1,\ldots,P\}$ and denote by $i$ the index of a generic nonleader car belonging to the $p$-th population. Denote the index of the first car ahead of the $i$-th one, belonging to population $p$, simply by $i+1$. Denote by $k$ the index of a generic car belonging to any population, including the $p$-th one. Denote the index of the car in front of the $k$-th car by $k+1$. Finally, denote by $k_i$ the index of a generic car belonging to any population between car $i$ (included) and car $i+1$ (excluded), see Fig.\ \ref{fig:proof}(top); and by $i_k$ the index of the first car belonging to population $p$ on the left of the car $k$ ($i_k=k$ if the $k$-th car belongs to population $p$), see Fig.\ \ref{fig:proof}(bottom). Assume for simplicity that the leftmost car belongs to population $p$ so that $i_k$ is well defined.

Denote the initial density for the vehicles of population $p$ by $\bar\mu\in R\cap BV(\R;[0,1])$ (dropping the superscript $p$ for simplicity), and denote the corresponding initial positions of vehicles by $\{\bar{y}_i\}_i$. Denote the positions of the $n_p$ vehicles of population $p$ at time $t$ by $\{y_i(t)\}_i$ and the positions of vehicles of any population ($p$-th included) at time $t$ by $\{y_k(t)\}_k$ (note that their number is variable in time). 
Finally let us denote by $T>0$ the final time for the simulation and consider a test function $\varphi\in C^\infty_0((-\infty,T]\times\R;\R)$. 

We are now ready to prove the correspondence between the microscopic follow-the-leader model (\ref{FTL_netw_paths}) and the macroscopic model (\ref{LWRlimit}), for any fixed $p\in\{1,\ldots,P\}$. 
Recalling \cite[Formula (4.5)]{bressanbook} the weak (integral) form of the equation (\ref{LWRlimit})  with the initial condition $\mu(0,\cdot)=\bar\mu$, we compute
\begin{equation*}
\begin{split}
& I_n:=\int_0^T\int_\R\big(\mu_n\varphi_t+\mu_n v^*(\omega_n)\varphi_x\big)dxdt+
\int_\R \bar\mu(x)\varphi(0,x)dx=\\
& \int_0^T\int_\R\mu_n\Big(\varphi_t+v^*(\omega_n)\varphi_x\Big)dxdt+
\int_\R \mu_n(0,x)\varphi(0,x)dx+\\ 
& \hskip8cm
\int_\R\big(\bar\mu(x)-\mu_n(0,x)\big)\varphi(0,x)dx.
\end{split}
\end{equation*}
By definition, the densities $\mu_n$ and $\omega_n$ are constant in the interval $[\ykt,\ykkt)$ for any $t$. Therefore we have
\begin{equation*}
\begin{split}
I_n=\int_0^T\sum_k\frac{\ell_n}{\dikt}\int_{\ykt}^{\ykkt}\Bigg(\varphi_t+\underbrace{v^*\left(\frac{\ell_n}{\Dkt}\right)}_{=\dot y_k(t)}\varphi_x\Bigg)dxdt+
\sum_i\frac{\ell_n}{\bardi}\int_{\baryi}^{\baryii} \varphi(0,x)dx+\\
\int_\R\big(\bar\mu(x)-\mu_n(0,x)\big)\varphi(0,x)dx,
\end{split}
\end{equation*}
where $\bardi=\bar y_{i+1}-\bar y_i$. 
Denoting by $\dot\varphi$ the total derivative of $\varphi$ w.r.t.\ time, we have, for $x\in[\ykt,\ykkt)$,
\begin{equation*}
\begin{split}
& |\varphi_t(t,x)+\dykt\varphi_x(t,x)-\dot\varphi(t,\ykt)|=\\
& |\varphi_t(t,x)+\dykt\varphi_x(t,x)-\varphi_t(t,\ykt)-\dykt\varphi_x(t,\ykt)|\leq \\
& \Cphi(1+\vmax)|\ykkt-\ykt|
\end{split}
\end{equation*}
where the quantity $\Cphi:=\|\varphi\|_{C^2}$ uniformly bounds from above the modulus of $\varphi$ and all its derivatives up to second order.
Equivalently, we can write
$$
\varphi_t(t,x)+\dykt\varphi_x(t,x)=\dot\varphi(t,\ykt))+\Ouno(\ykkt-\ykt).
$$
Using the latter estimate in the expression for $I_n$ above and defining
$$
A:=\sum_i\frac{\ell_n}{\bardi}\int_{\baryi}^{\baryii} \varphi(0,x)dx,\qquad
B:=\int_\R\big(\bar\mu(x)-\mu_n(0,x)\big)\varphi(0,x)dx,
$$
we get
\begin{equation*}
\begin{split}
I_n-A-B=\int_0^T\sum_k\frac{\ell_n}{\dikt}\int_{\ykt}^{\ykkt}\dot\varphi(t,\ykt))dxdt+\hskip4.5cm\\
\int_0^T\sum_k\frac{\ell_n}{\dikt}\int_{\ykt}^{\ykkt}\Ouno(\ykkt-\ykt)dxdt=
\end{split}
\end{equation*}
\begin{equation*}
\begin{split}
& \ell_n\int_0^T\sum_k\dot\varphi(t,\ykt)\frac{\Dkt}{\dikt}dt+
\Ouno\ell_n\int_0^T\sum_k\underbrace{\frac{\Dkt}{\dikt}}_{\leq 1}\Dkt dt=\\
& \ell_n\int_0^T\sum_i\sum_{k_i\in[i,i+1)}\dot\varphi(t,\ykit)\frac{\Dkit}{\dit}dt+
\Ouno\ell_n\int_0^T\sum_k\Dkt dt=
\end{split}
\end{equation*}
\begin{equation*}
\begin{split}
\ell_n\int_0^T\sum_i\sum_{k_i\in[i,i+1)}\Big(\dot\varphi(t,\yit)+\Ouno(\yiit-\yit)\Big)\frac{\Dkit}{\dit}dt+\hskip3.2cm\\
\Ouno\ell_n\int_0^T\sum_i\sum_{k_i\in[i,i+1)}\Dkit dt=
\end{split}
\end{equation*}
\begin{equation*}
\begin{split}
\ell_n\sum_i\int_0^T\Big(\dot\varphi(t,\yit)+\Ouno(\yiit-\yit)\Big)\underbrace{\sum_{k_i\in[i,i+1)}\frac{\Dkit}{\dit}}_{=1}dt+\hskip3.1cm\\
\Ouno\ell_n\sum_i\int_0^T\sum_{k_i\in[i,i+1)}\Dkit dt=
\end{split}
\end{equation*}
\begin{equation*}
\begin{split}
& \ell_n\sum_i\int_0^T\dot\varphi(t,\yit)dt+
\Ouno\ell_n\sum_i\int_0^T\dit dt+
\Ouno\ell_n\sum_i\int_0^T\dit dt=\\
& \ell_n\sum_i\int_0^T\dot\varphi(t,\yit)dt+
\Ouno\ell_n\sum_i\int_0^T\dit dt.
\end{split}
\end{equation*}
At this point we got rid of the populations $q$'s, $q\neq p$, and thus we are back to the case of a single road with a single population of vehicles. By the way, we have also resolved the discontinuity issues around the junctions. Therefore, from now on the proof follows the one in \cite{colombo2014RSMUP}. Recalling that $\varphi(T,\cdot)=0$, we have 
\begin{equation*}
\begin{split}
& I_n=-\ell_n\sum_i\varphi(0,\baryi)+\Ouno\ell_n\int_0^T (y_{n_p}(t)-y_1(t))dt+A+B=\hskip3cm
\end{split}
\end{equation*}
\begin{equation*}
\begin{split}
-\sum_i\frac{\ell_n}{\bardi}\int_{\baryi}^{\baryii} \varphi(0,\baryi)dx+
\Ouno\ell_n\int_0^T (y_{n_p}(t)-y_1(t))dt+\hskip3.5cm\\
\sum_i\frac{\ell_n}{\bardi}\int_{\baryi}^{\baryii} \varphi(0,x)dx+B=\\
\sum_i\frac{\ell_n}{\bardi}\int_{\baryi}^{\baryii} \Big(\varphi(0,x)-\varphi(0,\baryi)\Big)dx+
\Ouno\ell_n\int_0^T (y_{n_p}(t)-y_1(t))dt+B.\hskip1.8cm
\end{split}
\end{equation*}
To proceed further, let us define $\Lambda:=\bar y_{n_p}-\bar y_1$ and note that it does not depend on $n_p$ or $n$. We have
\begin{equation*}
\begin{split}
I_n=\sum_i\frac{\ell_n}{\bardi}\Ouno\int_{\baryi}^{\baryii}\bardi dx+
\Ouno\ell_n(\Lambda+\vmax T)T+\hskip4.4cm\\
\int_\R\big(\bar\mu(x)-\mu_n(0,x)\big)\varphi(0,x)dx=\\
\Ouno\ell_n\Lambda+
\Ouno\ell_n(\Lambda+\vmax T)T+
\int_\R\big(\bar\mu(x)-\mu_n(0,x)\big)\varphi(0,x)dx.\hskip2.5cm
\end{split}
\end{equation*}
Since $\mu_n(0,\cdot)\to\bar\mu$ (see \cite[Prop.\ 2.1]{colombo2014RSMUP}) and $\ell_n\to 0$, all the terms in the latter quantity vanish as $n\to+\infty$. 
Moreover there exists $\mu$ 
such that $\lim_{n\to\infty} \mu_n(t,x)=\mu(t,x)$ almost everywhere (see \cite[Th.\ 3]{difrancesco2015ARMA}).
In conclusion, we have proved the following 
\begin{theorem}\label{teo:mainresult}
Let (\ref{link_v_w}) and (\ref{properties_of_v}) hold, and consider the extensions (\ref{def:w*}) and (\ref{def:v*}). Fix $p\in\{1,\ldots,P\}$, and choose $\bar\mu^p=\mu^p(0,\cdot)\in R~\cap~BV(\R;[0,1])$ and $\bar{\mathbf{y}}^p=E_n[\bar\mu^p]$. 
Let $\mathbf{y}^p(\cdot)$ be the solution of (\ref{FTL_netw_paths}) with initial condition $\bar{\mathbf{y}}^p$ (assuming the other solutions $\mathbf{y}^q(\cdot)$, $q \neq p$ are given).
Define $\mu_n^p(t,\cdot)=C_n[\mathbf{y}^p]$. Then, $\mu^p(t,x):=\lim_{n\to +\infty}\mu_n^p(t,x)$ is a weak solution to (\ref{LWRlimit}) (assuming the other solutions $\mu^q$, $q\neq p$ are given) with initial condition $\bar\mu^p$.
\end{theorem}

Now we note that the system of PDEs (\ref{LWRlimit}) is nothing else that the multi-path model studied in \cite{bretti2014DCDS-S, briani2014NHM}, 
\begin{equation}\label{MPstandard}
\partial_t\mu^p+\partial_x(\mu^pv(\omega^p))=0, \qquad p=1,\ldots,P,\quad t>0,\quad x\in\text{path } p
\end{equation} 
where $\omega^p(t,x)$ is the sum of all densities $\mu^q(t,x)$, $q=1,\ldots,P$, living at time $t$ at the point $x$ along path $p$ and $v$ is the velocity of the cars.  
It is a system of $P$ conservation laws with discontinuous flux which provides an alternative method to deal with traffic flow on networks. It is characterized by the fact that junctions are embedded in the equations themselves, so that the dynamics at junctions have not to be resolved separately by \emph{ad hoc} procedures, like, e.g., the maximization of the flux. It appears that the model maximizes the flux but under the additional constraints of equidistribution of the fluxes from the incoming roads (i.e.\ the incoming roads have assigned the same priority). 

\begin{remark}
We proved the correspondence in the limit between the microscopic and the macroscopic model \emph{population by population}, i.e.\ fixing $p$ and assuming the other solutions to be given. This is different from considering the convergence of the fully coupled system of ODEs (\ref{FTL_netw_paths}) to the system of PDEs (\ref{LWRlimit}). On the other hand, a complete study is much harder because the microscopic vector field is not regular and because theoretical results for the multi-path model are still missing.
\end{remark}
%
%
%
%
\section{Numerical approximation}\label{sec:numerics}
In this section we give some details about the discretization of the microscopic and macroscopic equations.

\subsection{The follow-the-leader model} For numerical purposes, it is convenient to discretize the follow-the-leader model as it is described in section \ref{sec:FTLnetw_natural}, rather than its equivalent formulation on paths. 

First we set $\wmax=1$ and 
\begin{equation}\label{greeshields}
v(\rho)=1-\rho, \qquad w(\delta)=1-\frac{\ell_n}{\delta}.
\end{equation}
Then, we introduce a time step $\Delta t>0$ and we denote by $y_i^s$, $s=0,1,2,\dots$ the approximate position $y$ of the vehicle $i$ at time $t^s=s\Delta t$.
The discretization is obtained by the explicit Euler scheme 
\begin{equation}\label{FTLdiscrete}
y^{s+1}_i=y^s_i+\Delta t \ w^*(\hat\delta_i(t^s)), \qquad s=0,1,2,\dots
\end{equation}
for any car $i=1,\ldots,n$, and admissible initial condition $y_i^0=\bar y_i$.

When a vehicle crosses the junction and goes beyond it, let say $y_i(t^s)=L^{\RR(i,t^s)}+\epsilon$, it is assigned to the next road on the basis of its preference, and its position is updated as $y_i(t^s)=\epsilon$. We assume that no cars enter the road network after the initial time and that cars leave the network once they have reached a destination node. 

\medskip 

At the discrete level, a CFL-like condition is needed to ensure that vehicles do not bump into the one in front. Far from the junction, this condition is given by
\begin{equation}\label{CFL_FTL_gen}
\Delta t \ w(\delta_i)<\delta_i \qquad \forall i.
\end{equation}
Substituting (\ref{greeshields}) in (\ref{CFL_FTL_gen}) we get
\begin{equation}\label{CFL_FTL_partic}
\Delta t = \Delta t(n) < \min_i \frac{\delta_i^2}{\delta_i-\ell_n}.
\end{equation}
Unfortunately, this condition does not prevent vehicle overlapping near the junction, as already discussed in section \ref{sec:FTLnetw_natural}. In a discrete setting, the overlapping can happen also \emph{after} the junction, more precisely it can happen either in the space interval $(L^{\RR}-\ell_n,L^{\RR}]$ of one incoming road $\RR$ or in the space interval $[0,\Delta t \ \wmax)$ of one outgoing road. $\Dt \ \wmax$ is the maximal distance one vehicle can travel in one time step. It can be seen that in the discrete setting the overlapping zone shrinks to $\emptyset$ for $n\to +\infty$ and $\Delta t \to 0$, and not only for $n\to +\infty$. 
However, the overlapping is still ``mild'', in the sense that the number of bumped vehicles is bounded by the number of incoming roads, see section \ref{sec:FTLnetw_natural}. 

\medskip

\subsection{The multi-path model} For the numerical discretization of the multi-path system (\ref{LWRlimit}) we employ the Godunov-based scheme introduced in \cite{bretti2014DCDS-S}. For the sake of completeness, we briefly recall the scheme. 
Along each path $p$, a space discretization is considered, defining a space step $\Delta x>0$ and space nodes $x_k:=k\Delta x$, $k=0,\ldots,N_x^p$, where $N_x^p$ is the number of space nodes on path $p$. We denote by $\mu_k^{s,p}$ and $\omega_k^{s,p}$ the approximations of the values $\mu^p(t^s,x_k)$ and $\omega^p(t^s,x_k)$, respectively.

The scheme reads, for any $k=0,\ldots,N_x^p$ and $s=0,1,2,\ldots$,
\begin{equation}\label{scheme_multipath}
\mu_{k}^{s+1,p} =\mu_{k}^{s,p}-
\frac{\Delta t}{\Delta x}\left(\frac{\mu_{k}^{s,p}}{\omega_{k}^{s,p}} \ g(\omega_{k}^{s,p},\omega_{k+1}^{s,p})-
\frac{\mu_{k-1}^{s,p}}{\omega_{k-1}^{s,p}}\ g(\omega_{k-1}^{s,p},\omega_{k}^{s,p})
\right)
\end{equation}
with initial conditions $\mu^{0,p}_k=\bar\mu^p(x_k)$, and
$$
\omega^{s,p}_k=\sum_{q=1}^{P}\mu^{s,q}_k.
$$
Note that $\omega^{s,p}_k$ is the sum of all densities living at $x_k$ (along path $p$) at time $t^s$.
The function $g$ is instead the classical Godunov flux
$$
g(\rho_-,\rho_+):=
\left\{
\begin{array}{ll}
\min\{f(\rho_-),f(\rho_+)\}, & \textrm{if } \rho_-\leq \rho_+ \\
f(\rho_-), & \textrm{if } \rho_->\rho_+ \textrm{ and } \rho_-<\sigma \\
f(\sigma), & \textrm{if } \rho_->\rho_+ \textrm{ and } \rho_- \geq \sigma \geq \rho_+ \\
f(\rho_+), & \textrm{if } \rho_->\rho_+ \textrm{ and } \rho_+>\sigma
\end{array}
\right. 
$$
where $\sigma:=\arg\max_{\rho\in[0,1]} f(\rho)$. Recalling (\ref{greeshields}), we have $f(\rho)=\rho(1-\rho)$ and $\sigma=\frac12$.


We employ Dirichlet-zero boundary conditions $\mu_k=0$ at origin and destination nodes.

\begin{remark} 
The multi-path model, at least at the discrete level, does not allow the total density to be larger than 1 \cite[Sect.\ 3.2]{briani2014NHM}. 
This seems in contradiction with the possibility of overlapping in the microscopic model. Actually this is not, because the overlapping is ``mild'', i.e.\ it concerns at most a finite number of cars, and then it is invisible at macroscopic level. 
However, a counterpart of the overlapping features at macroscopic level exists: the cells after the junctions act as a sort of \emph{buffer} \cite[Sect.\ 5.2]{briani2014NHM}, which is able to gather a flux of vehicles larger than maximal one, which is $\max_\rho \{\rho v(\rho)\}$. 
The actual maximal flux is instead $N_{\textup{inc}}\max_\rho \{\rho v(\rho)\}$, where $N_{\textup{inc}}$ is the number of incoming roads.
\end{remark}

\section{Numerical tests}\label{sec:tests}
\begin{figure}[h!]
\begin{center}
\begin{overpic}
[width=0.49\textwidth]{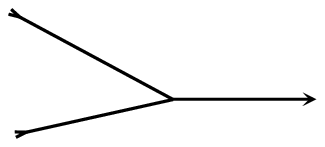}
\put(20,38){Road 1} \put(18,15){Road 2} \put(67,20){Road 3}
\end{overpic}
\begin{overpic}[width=0.49\textwidth,tics=10]{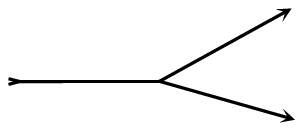}
\put(20,22){Road 1} \put(60,35){Road 3} \put(70,15){Road 4}
\end{overpic}
\\
\begin{overpic}[width=0.49\textwidth,tics=10]{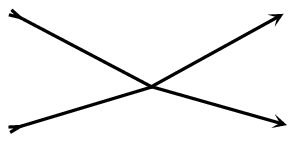}
\put(20,38){Road 1} \put(14,16){Road 2} \put(58,35){Road 3} \put(70,17){Road 4}
\end{overpic}
\end{center}
\caption{Top-left: a merge (3 roads and 2 paths). Top-right: a diverge (3 roads and 2 paths). Bottom: a junction with 2 incoming and 2 outgoing roads (4 roads and 4 paths).}
\label{fig:networks_for_tests}
\end{figure}

In this section we study the micro-to-macro limit by means of some numerical tests, finding a good agreement with theoretical findings of the previous sections. We consider the case of a simple network with 1 junction and (i) a merge, (ii) a diverge, and (iii) 2 incoming roads and 2 outgoing roads, see Fig.\ \ref{fig:networks_for_tests}.

In the cases (ii) and (iii), additional parameters are needed to describe the behavior of the vehicles at the junction. They are usually referred to as \emph{distribution coefficients} and specify the percentage of vehicles which wish to turn to the left and to the right. 
Since we are considering here simple networks with only 1 junction, specifying the car's turning conduct at junction is equivalent to define the car's path on the network. In the case of more complex networks with more than 1 junction, we should instead assign the whole sequence of decisions each car makes at junctions. If we stick with per-junction coefficients, thus losing the global behavior of drivers, we expect convergence, in the limit, to the hybrid version of the multi-path model described in \cite[Remark 1]{briani2014NHM}.

Let us denote the distribution coefficients by
$
\mathcal{P}_{a\to b}
$, 
where $a$ is the index of the incoming road and $b$ is the index of the outgoing road. 
Clearly it is required that 
$$
\sum_b \mathcal{P}_{a\to b}=1\qquad \forall a.
$$
In the follow-the-leader model, each vehicle traveling along road $a$ is assigned to the road $b$ with probability $\mathcal P_{a\to b}$. In the multi-path model, these coefficients are used to define properly the partial densities $\mu^p$'s at the initial time given the initial condition for the total density. 

To improve readability and ease the comparison, we only show the total density, and we redefine it on single roads, rather than on paths. We also come back to the original notation $\rho$, denoting by $\rho^{\RR}$ the total density on each road $\RR=1,\ldots,N_\RR$. We assume that all roads have the same length $L^\RR=4\times 10^3$, and we divide them in $N_x^{\textsc{\RR}}=100$ space nodes. 

Comparison with the solution of the multi-path scheme is obtained by computing the average density of the microscopic vehicles, defined for each space node $k$ and time step $s$ as
\begin{equation}\label{lagrangian_density}
\Psi_k^s[\mathbf{y}^s]:=
\frac{\ell_n}{\Delta x}\#{\{i~:~y_i^s\in [x_k,x_{k+1})\}},
\end{equation}
where $\#A$ denotes the cardinality of any set $A$. This quantity corresponds to the total length of cars found in each space cell divided by the length of the cell. 

\subsection{Merge}
\begin{figure}[t!]
\begin{center}
\includegraphics[width=0.46\textwidth]{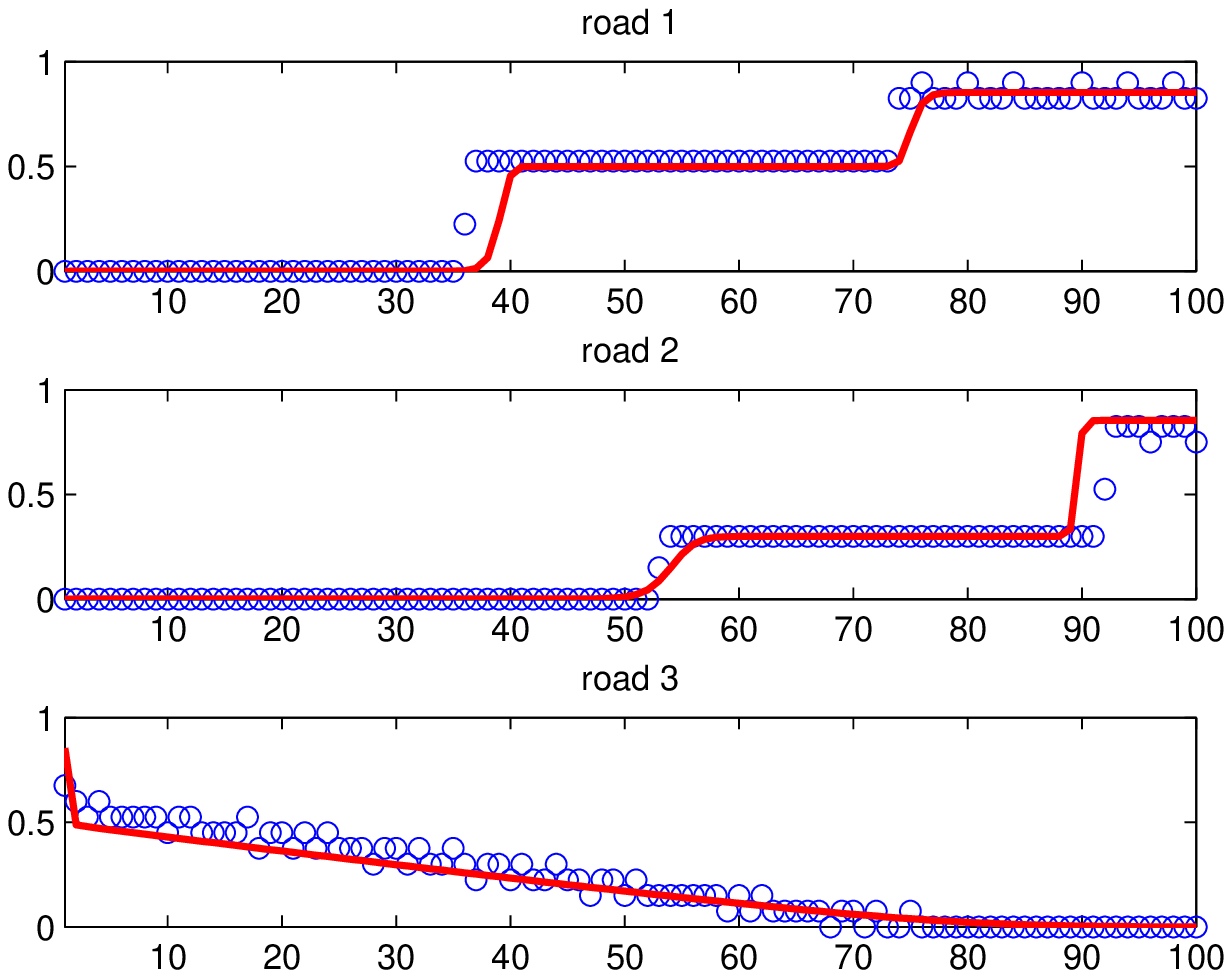} \qquad
\includegraphics[width=0.46\textwidth]{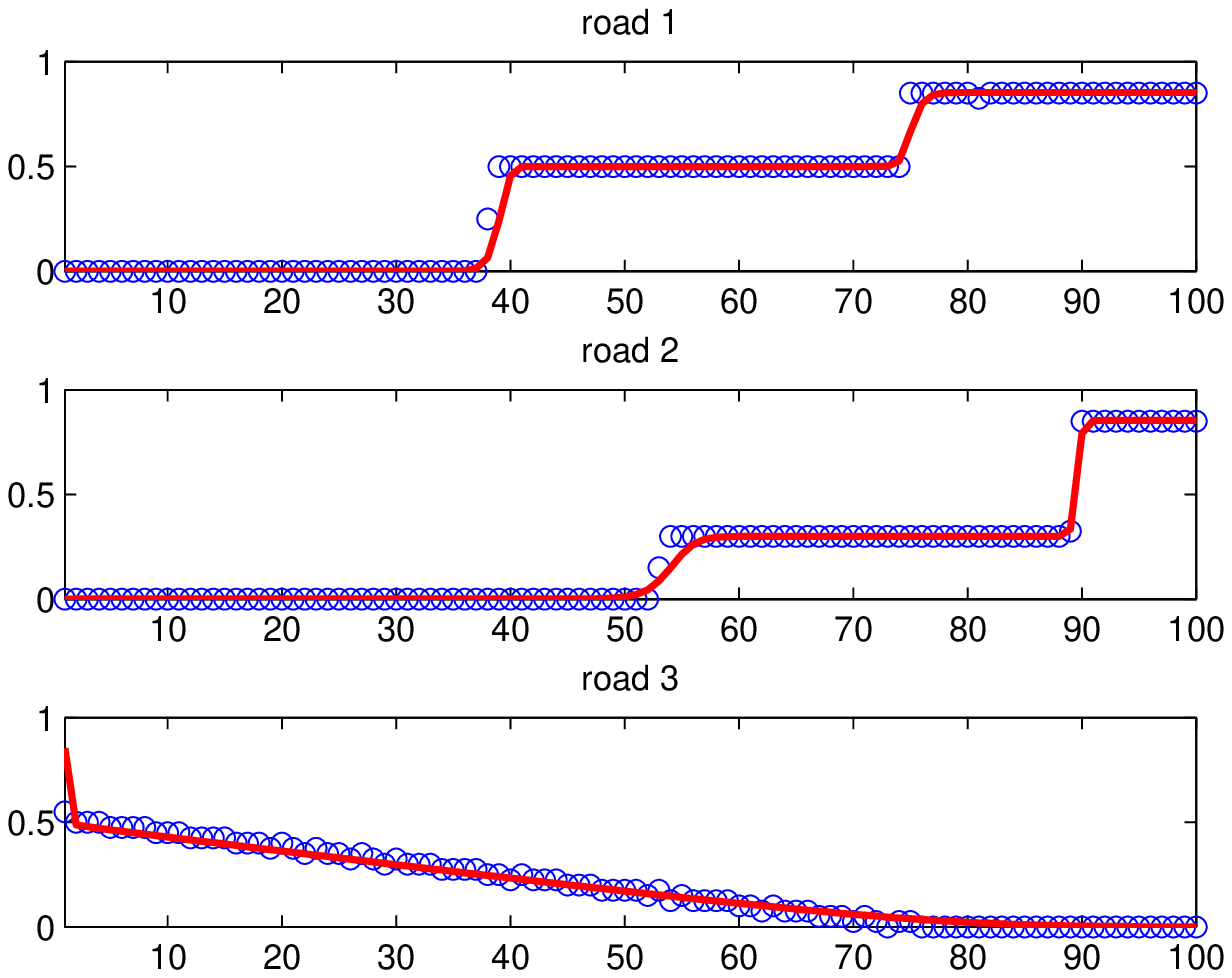}
\end{center}
\caption{Merge, result of the simulation at final time. Total macroscopic density redefined on roads (red line) and density of microscopic vehicles $\Psi$ (blue circles). Left: $\ell_n=3$ (max number of cars per cell $= 13.\bar 3$), $\Delta t=3$. Right: $\ell_n=1$ (max number of cars per cell $= 40$), $\Delta t=0.2$.}
\label{fig:2in1_results}
\end{figure}

In this section we consider a network with three roads and one junction, with two incoming roads and one outgoing road. We denote by $\rho^1$, $\rho^2$ and $\rho^3$ the density on the first incoming road, the second incoming road and the outgoing road, respectively (Fig.\ \ref{fig:networks_for_tests}). 
Initial conditions are
$$
\rho^1(0,x)=0.5, \qquad \rho^2(0,x)=0.3, \qquad \rho^3(0,x)=0 \qquad \forall x.
$$
We run the simulation until the final time $T=3\times 10^3$. Fig.\ \ref{fig:2in1_results} shows the density computed by the multi-path model and the follow-the-leader model for two different choices of $\ell_n$ and $\Delta t$.
Since the outgoing road is not able to gather the fluxes coming from the incoming roads, two queues are formed and propagate back along the incoming roads (with different speed). Incoming fluxes are equidistributed and queues have the same level of density, equal to $\rho^*$ such that $\rho^* v(\rho^*)=\frac12 \sigma v(\sigma)$.
Since we compare two numerically approximate densities, we expect a perfect match only for $n\to +\infty$ and $\Delta t, \Delta x\to 0$. That said, numerical evidence confirms the convergence results, and also shows that the microscopic scheme is not diffusive as instead it is the macroscopic scheme. This is perfectly visible across the discontinuities.
\subsection{Diverge}
\begin{figure}[t!]
\begin{center}
\includegraphics[width=0.46\textwidth]{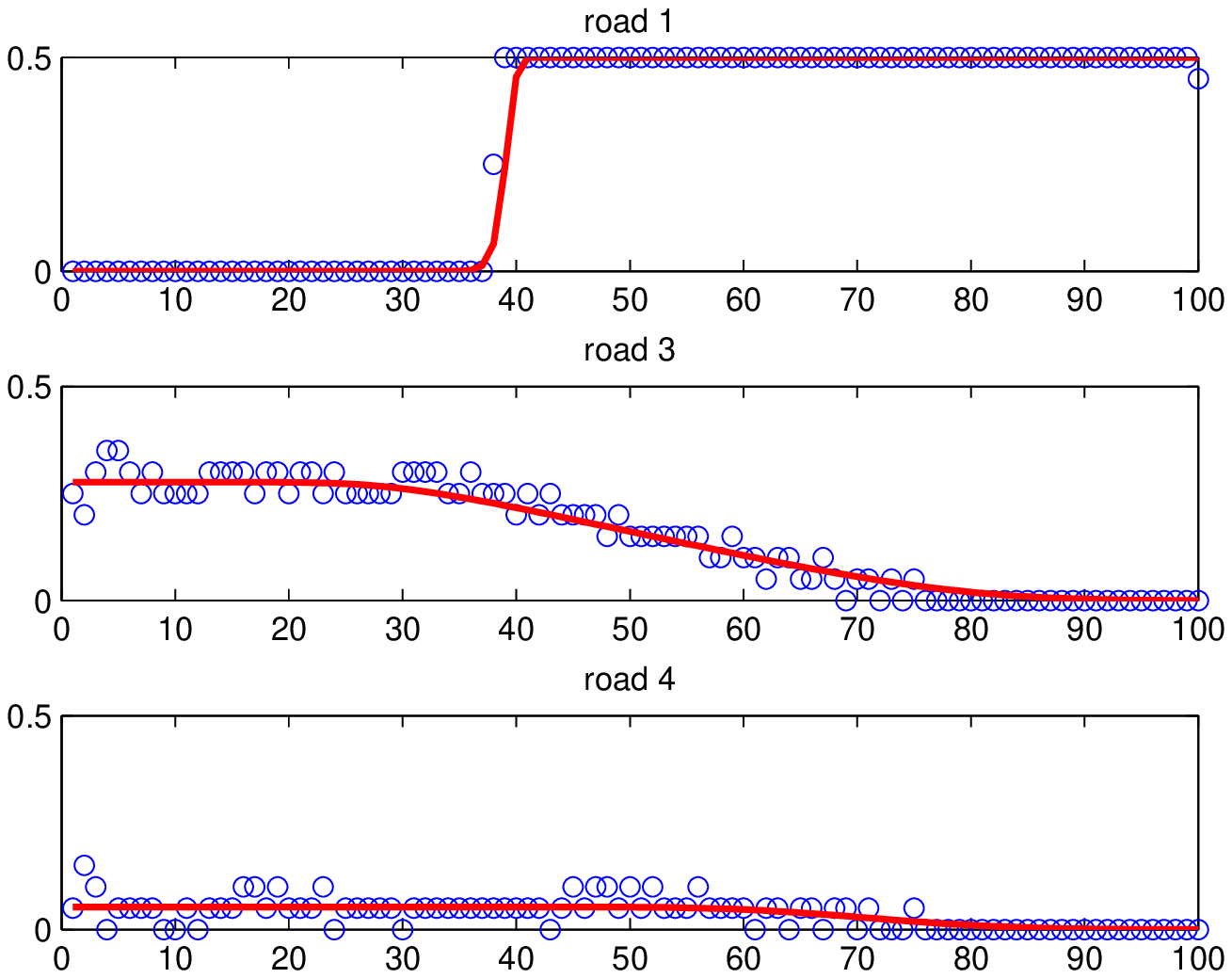}\qquad
\includegraphics[width=0.46\textwidth]{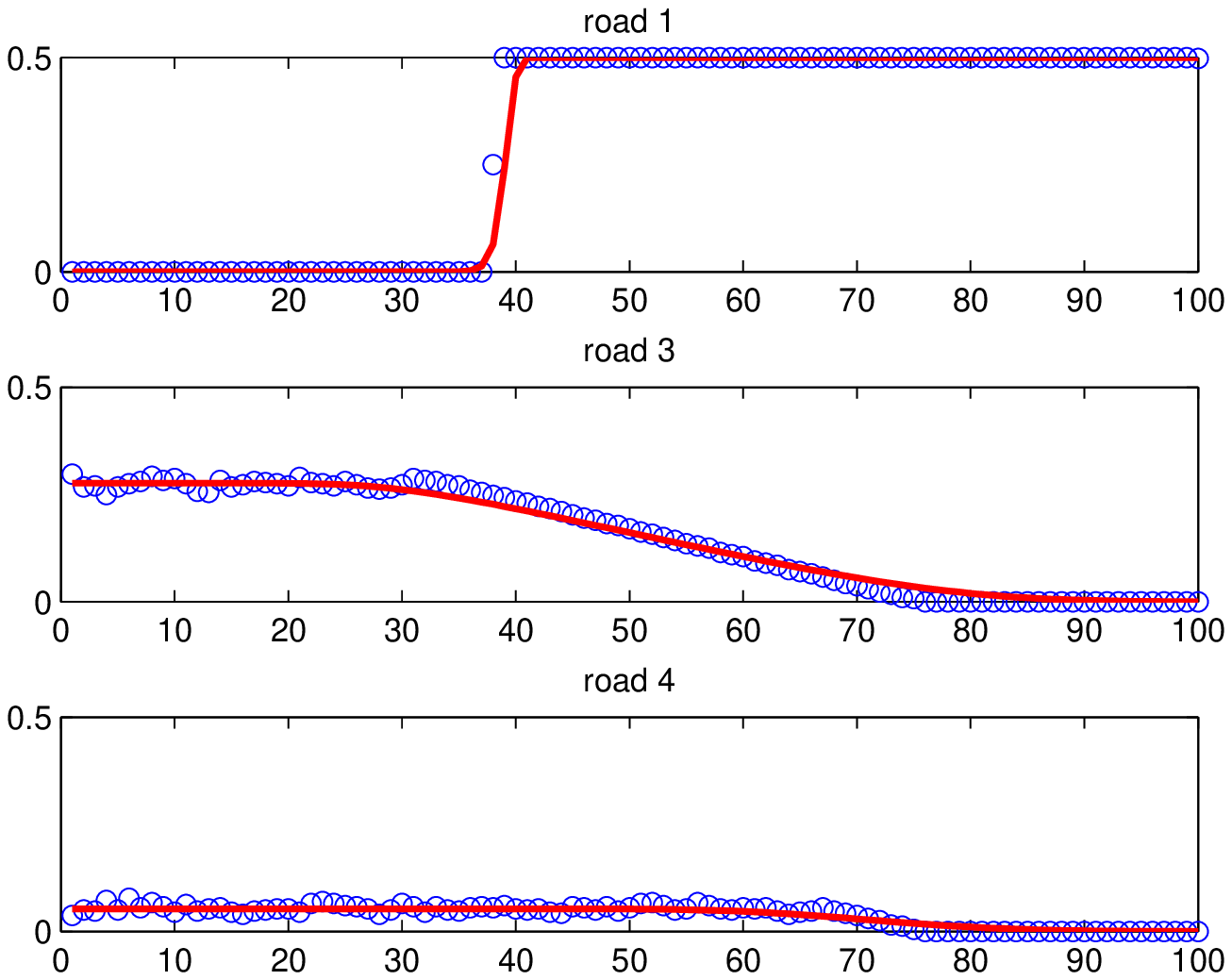}
\end{center}
\caption{Diverge, result of the simulation at final time. Total macroscopic density redefined on roads (red line) and density of microscopic vehicles $\Psi$ (blue circles). Left: $\ell_n=2$ (max number of cars per cell $= 20$), $\Delta t=4$. Right: $\ell_n=0.1$ (max number of cars per cell $= 400$), $\Delta t=0.25$.}
\label{fig:1in2_results}
\end{figure}

In this section we consider a network with three roads and one junction, with one incoming road and two outgoing roads. We denote by $\rho^1$, $\rho^3$ and $\rho^4$ the density on the incoming road, the first outgoing road and the second outgoing road, respectively (Fig.\ \ref{fig:networks_for_tests}). 
Initial conditions are
$$
\rho^1(0,x)=0.5, \qquad \rho^3(0,x)=0, \qquad \rho^4(0,x)=0 \qquad \forall x,
$$
and distribution coefficients are
$$
\mathcal{P}_{1\to 3}=0.8, \qquad \mathcal{P}_{1\to 4}=0.2.
$$
We run the simulation until the final time $T=3\times 10^3$. Fig.\ \ref{fig:1in2_results} shows the density computed by the multi-path model and the follow-the-leader model for two different choices of $\ell_n$ and $\Delta t$.
As expected, the density splits among the two outgoing roads. The corresponding flux splits with ratio $\frac{0.8}{0.2}$ and mass is conserved. Again, numerical evidence confirms the convergence results, however, the convergence is much slower than the previous case. This is probably due to the presence of the distribution coefficients $\mathcal{P}$ which introduce stochasticity in the system. Note that we are showing the outcome of a single run and not the average of many runs.
\subsection{Junction with two incoming and two outgoing roads}
\begin{figure}[h!]
\begin{center}
\includegraphics[width=0.7\textwidth]{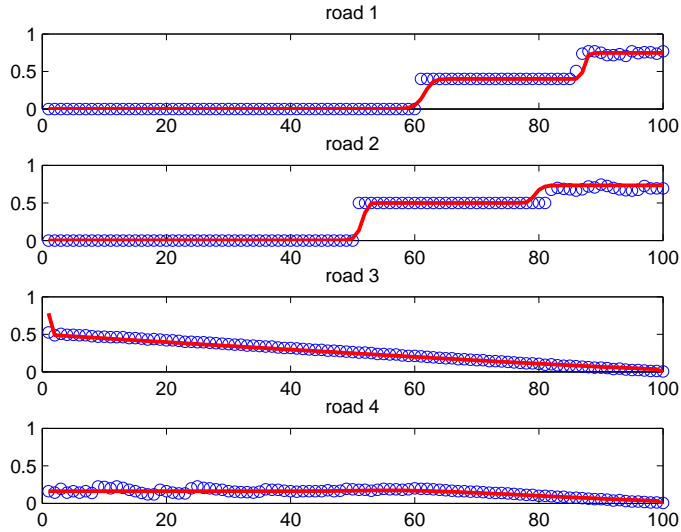}
\end{center}
\caption{2-in-2 junction, result of the simulation at final time. Total macroscopic density redefined on roads (red line) and density of microscopic vehicles $\Psi$ (blue circles). $\ell_n=0.25$ (max number of cars per cell $= 160$), $\Delta t=0.1$.}
\label{fig:2in2_results}
\end{figure}

In this section we consider a network with four roads and one junction, with two incoming roads and two outgoing roads. We denote by $\rho^1$, $\rho^2$, $\rho^3$ and $\rho^4$ the density on the first incoming road, the second incoming road, the first outgoing road, and the second outgoing road, respectively (Fig.\ \ref{fig:networks_for_tests}). 
Initial conditions are
$$
\rho^1(0,x)=0.4, \qquad \rho^2(0,x)=0.5, \qquad \rho^3(0,x)=0, \qquad \rho^4(0,x)=0 \qquad \forall x,
$$
and distribution coefficients are
$$
\mathcal{P}_{1\to 3}=0.7, \qquad \mathcal{P}_{1\to 4}=0.3, \qquad \mathcal{P}_{2\to 3}=0.6, \qquad \mathcal{P}_{2\to 4}=0.4.
$$
We run the simulation until the final time $T=4\times 10^3$. Fig.\ \ref{fig:2in2_results} shows the density computed by the multi-path model and the follow-the-leader model.
Again, numerical evidence confirms the convergence results, however, the convergence is even slower than the previous case.
%
%
%
%
%
\section*{Conclusions and future work}
In this paper we have investigated the many-particle limit for a natural extension of the follow-the-leader model on networks. Numerical tests have shown a slow convergence to the limit solution. This means that a large number of particles is required to get a precision comparable with that of the associate macroscopic model. In addition, the vector storing cars' indices must be reordered after every time step in order to find the car in front of any car, so that a long computational time is needed.  
In conclusion, we discourage the use of the microsimulator but our results justify any multiscale approaches where the micro and the macro scale live together and exchange information.

The convergence result proven in this paper is only one of the possible micro-to-macro limits which can be investigated in the framework of traffic flow models on networks. 
First of all, one could define the proper follow-the-leader model which corresponds to, e.g., the LWR model with maximization of flux at junctions. This requires a nontrivial management of the junctions, in which an authority is able to decide who passes the junction and when. 

Second, one can consider analogous techniques for second-order models. In this case the problem is completely open because a second-order multi-path does not yet exist and then it can be difficult to find the right limit model.

Finally, it could be interesting the investigation of the many-particle limit for meso-to-macro models. To this end, useful references could be the papers \cite{fermo2015M3AS}, which deals with kinetic models on networks, and \cite{fermo2014DCDS-S} about fundamental diagrams for kinetic equations.

\section*{Acknowledgments}
The authors wish to thank F. S. Priuli, M. Falcone, and B. Piccoli for the useful discussions and suggestions.


\bigskip
\bigskip

\end{document}